[1994/12/01]
\documentclass{aomart}

\usepackage[english]{babel}
\usepackage{slashbox}
\usepackage{xcolor}
\usepackage{colortbl}
\definecolor{lgrey}{rgb}{0.8,0.8,0.8}

\theoremstyle{definition}
\newtheorem{defn}{Definition}[section]
\newtheorem{ass}{Assumption}

\title[About the proof of the Collatz conjecture]{About the proof of the Collatz conjecture}
\author{Carolin Z\"obelein}
\address{Friedrich-Alexander Universit\"at Erlangen-N\"urnberg, Germany \\
 Student of Physics Departement}
\email{carolin.zoebelein@physik.stud.uni-erlangen.de}
\givenname{Carolin}
\surname{Z\"obelein}

\begin{document}

\begin{abstract}
	I want to show one possibility to proof the Collatz conjecture, also called $3n+1$ conjecture, for any natural number $N \in \mathbb{N}$. For this, I limit my analysis on the direct odd follower of every natural odd number and show the connections between the already by one reached numbers and their followers, to find an recurrence over all ranges $[1,N_{i}]$, $i \in \mathbb{N}$ to proof the conjecture.
\end{abstract}

\maketitle
\tableofcontents

\section{Introduction}
\label{s:int}
The Collatz conjecture is given by the following number sequence for any natural number $n$
\begin{equation} n_{i+1} = \left\{
	\begin{array}{l@{\quad \quad}l}
	\frac{n_{i}}{2} & \textnormal{if $n_{i}$ is even} \\
	3n_{i}+1 & \textnormal{if $n_{i}$ is odd}
	\end{array}
\right. \label{eq:collatz} \end{equation} 
and the observation that it seems to be, that $\forall n \in \mathbb{N}$, this number sequence always ends in the only closed chain $4 \rightarrow 2 \rightarrow 1$.

To proof this observation it has to show, that, on the one hand, no other closed chain exists, and on the other hand, that it exists no open chain with endless length. 

This two characteristics of Collatz sequence will now shown in the following sections.

This paper bases upon ideas of the papers \cite{Work2s3s}, \cite{OnCp}, \cite{BouCc} and \cite{2011arXiv1108.4056L}.
\section{Number generation by recurrence of open chains}
\label{s:genopch}
\subsection{Product factors between sequenced numbers}
\label{ss:prodfact}
To check the trueness of the conjecture, we will, at first, not consider the concrete numbers of the sequence, instead, we will consider the product factors $f$ between $n_{i}$ and $n_{i+1}$, $i \in \mathbb{N}$, and their qualities.
We look at first on the case, if $n_{i}$ is even. So we have
\[ n_{i+1} = \frac{n_{i}}{2} \]
\[ n_{i+1} = f_{e} n_{i} \]
\begin{equation}\label{eq:faceven}
	f_{e,i(i+1)} = \frac{1}{2} \Leftrightarrow f_{e,(i+1)i} = 2
\end{equation}
And for the case, if $n_{i}$ is odd, we have
\[ n_{i+1} = 3n_{i} + 1 \]
\[ n_{i+1} = f_{o} n_{i} \]
\begin{equation}\label{eq:facodd}
	f_{o,i(i+1)} = \frac{3n_{i}+1}{n_{i}} \Leftrightarrow f_{o,(i+1)i} = \frac{n_{i}}{3n_{i}+1}
\end{equation}
\subsection{Product consideration of open and closed chains}
\label{ss:clch}
Now we are interested in the evolution of the total factor of the number sequence after a certain number of steps. For this, we will differ the number sequences in two possible variations. At first, have a look at, so called, \textit{open chains}.

Let
\[ n_{1} \rightarrow n_{2} \rightarrow \cdots \rightarrow n_{l-1} \rightarrow n_{l} \]
be a Collatz sequence, started by a random number $n_{1}$ and ending after $l-1 \in \mathbb{N}$ steps on an number $n_{l}$, $n_{1},... ,n_{l} \in \mathbb{N}$. 

For the following considerations, we will expect that $n_{1}$ always is odd. This assumption will not limit the result of our proof, since we can directly see that every even number always ends, after a certain number of steps, on an odd number. Thus, we only have to show the trueness of the conjecture for all odd numbers, to show it for all even numbers too. 

\begin{defn}[Open chain] We call a number sequence an open chain, if applies $n_{i} \neq n_{j}, i \neq j$ for $l-1$ steps. The chain consists of $e$ even and $o$ odd steps. Step means in this context to go from element $n_{i}$ to element $n_{i+1}$. So remark, $e$ and $o$ are not the counts of even respectively odd numbers, since the last number $n_{l}$ makes no step.  
\end{defn}

For the total factor between $n_{1}$ and $n_{l}$ after $l-1$ steps for an open chain, it follows
\begin{equation}\label{eq:proop}
	2^{e} \prod_{j=1}^{o} \frac{n_{j}}{3n_{j}+1} = \frac{n_{1}}{n_{l}} \Leftrightarrow \frac{1}{2^{e}} \prod_{j=1}^{o} \frac{3n_{j}+1}{n_{j}} = \frac{n_{l}}{n_{1}} 
\end{equation}
Remark, $j$ is not the $j$'th element. Now, it is the $j$'th odd number inside the sequence.

At next, have a look at, so called, \textit{closed chains}. We will start with the same sequence from just, but with a little difference.

\begin{defn}[Closed chain] We call a number sequence a closed chain, if applies $n_{l} = n_{1}$ and $n_{i} \neq n_{j}, i,j \in [2,l-1]$ and $i \neq j$, for $l$ steps. The chain consists of $\bar{e}$ even and $\bar{o}$ odd steps. So remark, now $\bar{e}$ and $\bar{o}$ are the counts of even respectively odd numbers, since the last number $n_{l}$ makes one step.
\end{defn}

For the total factor between $n_{1}$ and, after $l$ steps, as well $n_{1}$ for an closed chain, it follows
\begin{equation}\label{eq:procl}
	2^{\bar{e}} \prod_{j=1}^{\bar{o}} \frac{n_{j}}{3n_{j}+1} = 1 \Leftrightarrow \frac{1}{2^{\bar{e}}} \prod_{j=1}^{\bar{o}} \frac{3n_{j}+1}{n_{j}} = 1
\end{equation}
\subsection{Recurrence of open chains}
\label{ss:it_opch}
For the next proof steps we will take the following consideration of the Collatz sequence. We will start at the last element of sequence $n_{l}$ and will show, that if is $n_{l}=1$, we can reach all odd numbers higher than one by inverse Collatz sequence. For this, we will use the just made definition of open chains. 
\begin{equation}\label{eq:opone}
	2^{e} \prod_{j=1}^{o} \frac{n_{j}}{3n_{j}+1} = n_{1}
\end{equation}
We start with an exception. If we have no other odd numbers, apart from one, in our sequence, it follows with $o=0$
\begin{equation}\label{eq:opo0}
	2^{e_{0}} = n_{1} 
\end{equation}
This is the only open chain which has only even steps. So it will be not a result of the proof of the next sections and have to be separately add on in due time. Now let be $o=1$
\[ 2^{e_{1}} \frac{n_{1}}{3n_{1}+1} = n_{1} \]
and solve for $n_{1}$
\begin{equation}\label{eq:opo1}
	n_{1} = \frac{2^{e_{1}}-1}{3}
\end{equation}
Let be $o=2$
\[ 2^{e_{2}} \frac{n_{1}}{3n_{1}+1} \frac{n_{2}}{3n_{2}+1} = n_{1} \]
This we can rewrite with (\ref{eq:opo1}) to 
\[ 2^{e_{2}-e_{1}} \frac{n_{1}}{3n_{1}+1} n_{2} = n_{1} \]
so it follows for the new $n_{1}$ of a number sequence 
\begin{equation}\label{eq:opo2}
	n_{1} = \frac{2^{e_{2}-e_{1}}n_{2}-1}{3}
\end{equation}
Let be $o=3$
\[ 2^{e_{3}} \frac{n_{1}}{3n_{1}+1} \frac{n_{2}}{3n_{2}+1} \frac{n_{3}}{3n_{3}+1} = n_{1} \]
\[ 2^{e_{3}-(e_{2}-e_{1})} \frac{n_{1}}{3n_{1}+1} n_{2} = n_{1} \]
\begin{equation}\label{eq:opo3}
	n_{1} = \frac{2^{e_{3}-(e_{2}-e_{1})}n_{2}-1}{3}
\end{equation}
and so on. It follows, in general
\begin{equation}\label{eq:opngen}
	n_{1} = \frac{2^{x}n_{2}-1}{3}
\end{equation}
with $x \in \mathbb{N}$.
\section{General study of recurrence formula}
\label{s:stitfor}
\subsection{Numbers of recurrence formula}
\label{ss:nuitfor}
In previous section we got a general formula for recurrence of open chains. Now we will make some basic studies about this formula. Let us start with a view of the concrete series of numbers. 
The first question is, when exist integer solutions of (\ref{eq:opngen})? For this, we have to partition the whole set of odd natural numbers in three disjoint subsets. The first subset is the set of all odd numbers, which are integer multiple numbers of three. This subset has no solutions, since
\begin{equation}\label{eq:ssthr}
	n_{1} = \frac{2^{x}3j - 1}{3} = 2^{x}j - \frac{1}{3}
\end{equation} 
$j \in \mathbb{N}$. This means in this context, that, this numbers are generated but itself no generated new numbers. Within a Collatz sequence they don't have  odd predecessors. The second subset is given by
\begin{equation}\label{eq:sspev}
	n_{2,epow} = 6i_{epow} + 1
\end{equation}
$i_{epow} \in \mathbb{N}$, and additional $n_{2,epow}=1$, which solves the recurrence formula for even powers $x_{epow}$. The first solving numbers are listed in table \ref{tab:nopoweven}. 
\begin{table}
\centering
\begin{tabular}{|c||c|c|c|c|c|c|c|c|c|}\hline
	\backslashbox{$n_{2}$\kern-1em}{\kern-1em $x$} & $2$ & $4$ & $6$ & $8$ & $10$ & $12$ & $14$ & $16$ & $18$ \\ 
	\hline \hline $1$ & \textbf{1} & \textbf{5} & \cellcolor{lgrey} $21$ & $85$ & $341$ & \cellcolor{lgrey} $1365$ & $5461$ & $21845$ & \cellcolor{lgrey} $87381$ \\
	\hline $7$ & \cellcolor{lgrey} \textbf{9} & $37$ & $149$ & \cellcolor{lgrey} $597$ & $2389$ & $9557$ & \cellcolor{lgrey} $38229$ & $152917$ & $611669$ \\
	\hline $13$ & \textbf{17} & \cellcolor{lgrey} $69$ & $277$ & $1109$ & \cellcolor{lgrey} $4437$ & $17749$ & $70997$ & \cellcolor{lgrey} $283989$ & $1135957$ \\
	\hline $19$ & $25$ & $101$ & \cellcolor{lgrey} $405$ & $1621$ & $6485$ & \cellcolor{lgrey} $25941$ & $103765$ & $415061$ & \cellcolor{lgrey} $1660245$ \\
	\hline 
\end{tabular}\caption{Subset $n_{2,epow}=6i_{epow}+1$. Grey marked: all numbers which don't generate other numbers. Bold marked, all numbers of assumption range $[1,19]$ (see section \ref{s:itofN}).}\label{tab:nopoweven}\end{table}
The third subset is given by
\begin{equation}\label{eq:ssmod}
	n_{2,opow} = 6i_{opow} - 1
\end{equation}
$i_{opow} \in \mathbb{N}$, which solves the recurrence formula for odd powers $x_{opow}$. The first solving numbers are listed in table \ref{tab:nopoodd}.
\begin{table}
\centering
\begin{tabular}{|c||c|c|c|c|c|c|c|c|c|}\hline
	\backslashbox{$n_{2}$\kern-1em}{\kern-1em $x$} & $1$ & $3$ & $5$ & $7$ & $9$ & $11$ & $13$ & $15$ & $17$ \\
	\hline \hline $5$ & \cellcolor{lgrey} \textbf{3} & \textbf{13} & $53$ & \cellcolor{lgrey} $213$ & $853$ & $3413$ & \cellcolor{lgrey} $13653$ & $54613$ & $218453$ \\
	\hline $11$ & \textbf{7} & $29$ & \cellcolor{lgrey} $117$ & $469$ & $1877$ & \cellcolor{lgrey} $7509$ & $30037$ & $120149$ & \cellcolor{lgrey} $480597$ \\
	\hline $17$ & \textbf{11} & \cellcolor{lgrey} $45$ & $181$ & $725$ & \cellcolor{lgrey} $2901$ & $11605$ & $46421$ & \cellcolor{lgrey} $185685$ & $742741$ \\
	\hline 
\end{tabular} \caption{Subset $n_{2,opow}=6i_{opow}-1$. Grey marked: all numbers which don't generate other numbers. Bold, marked all numbers of assumption range $[1,19]$ (see section \ref{s:itofN}).} \label{tab:nopoodd} \end{table}
With this subsets and the recurrence formula, one can get all start numbers $n_{1}$ of Collatz sequences, generated from their direct follower. But thus, we have one exception in our recurrence. The case $(n_{2}, n_{1})=(1,1)$ gives us the recurrence view on the known closed chain $4 \rightarrow 2 \rightarrow 1 \rightarrow 4$. Here, the only odd number $1$ in the sequence is the factor between $n_{1}$ and $n_{l}$, too. So it iterates on itself and we don't can it use for the following recurrence considerations. 

Finally we can easy show that no two cases exist that generate the same number, since
\[ \frac{2^{x_{1}}n_{2,1}-1}{3} = \frac{2^{x_{2}}n_{2,2}-1}{3} \]
\[ 2^{x_{1}} n_{2,1} = 2^{x_{2}} n_{2,2} \]
\[ 2^{x_{1}-x_{2}} = \frac{n_{2,2}}{n_{2,1}} \]
is only solved for $x_{1}=x_{2} \Rightarrow n_{2,1} = n_{2,2}$, since $n_{2,1}, n_{2,2}$ are always odd. This quality is important for the next proof steps.
\subsection{Power relation between two $n_{2}$ number rows}
\label{ss:powrel}
Now we want consider the power relation between two different number rows for the same $n_{1}$. For this, let be
\[ n_{1,i} = \frac{2^{x_{i}} n_{2,i} - 1}{3} \]
and
\[ n_{1,j} = \frac{2^{x_{j}} n_{2,j} - 1}{3} \]
$i,j \in \mathbb{N}: i \neq j$. And solve
\[ 0 = n_{1,i} - n_{1,j} = \frac{2^{x_{i}} n_{2,i} - 1}{3} - \frac{2^{x_{j}} n_{2,j} - 1}{3} \]
\[ 0 = 2^{x_{i}} n_{2,i} - 2^{x_{j}} n_{2,j} \]
\[ x_{j} = \log_{2}\left( \frac{n_{2,i}}{n_{2,j}} 2^{x_{i}} \right) \]
\begin{equation}\label{eq:powrel}
	x_{j} = x_{i} + \log_{2} \left( \frac{n_{2,i}}{n_{2,j}} \right)
\end{equation}
\subsection{Odd powers: qualities for fixed $N$}
\label{ss:oddqual}
Let $N$, $N \in \mathbb{N}$, be a fixed odd number with
\begin{equation}\label{eq:N}
	N = \frac{2^{x_{N}}-1}{3}
\end{equation}
$x_{N} = 2k_{N}$. For consideration of odd powers we get for (\ref{eq:powrel}) with
\[ n_{2,i} = 1, x_{i} = 2k_{N} \]
\[ n_{2,j} = 6i_{opow} - 1, x_{j} = 2k_{j} - 1 \]
$k_{j} \in \mathbb{N}$.
\begin{equation}\label{eq:oddpowk}
	k_{j} = \left\lfloor k_{N} + \frac{1}{2} \log_{2} \left( \frac{1}{6i_{opow}-1} \right) + \frac{1}{2} \right\rfloor
\end{equation}
At first be $k_{N} \in \mathbb{N}: k_{N} > 1$ and look for solutions for $i_{opow}$. Equation (\ref{eq:oddpowk}) has integer solutions for
\[ \frac{1}{2} \log_{2}\left( \frac{1}{6i_{opow}-1} \right) = -\frac{1}{2}(2f-1) \]
$f \in \mathbb{N}$.
\[ 2^{-2f+1} = \frac{1}{6i_{opow}-1} \]
\begin{equation}\label{eq:oddpowi}
	i_{opow} = \left\lfloor \frac{1}{6} \left( \frac{1}{2^{-2f+1}} + 1 \right) \right\rfloor = \frac{1}{6} \left( \frac{1}{2^{-2f+1}} + 1 \right) - \frac{1}{2}
\end{equation}
and the maximum $i_{opow,max}$ for $f=k_{N}$
\begin{equation}\label{eq:oddpowimax}
	i_{opow,max} = \left\lfloor \frac{1}{6} \left( \frac{1}{2^{-2k_{N}+1}} + 1 \right) \right\rfloor = \frac{1}{6} \left( \frac{1}{2^{-2k_{N}+1}} + 1 \right) - \frac{1}{2}
\end{equation}
At second be $N=2p_{N}-1, p_{N} \in \mathbb{N}: p_{N} > 1$, any odd natural number. So
\begin{equation}\label{eq:anyN}
	k_{N} = \frac{1}{2} \log_{2} \left( 3N+1 \right) = \frac{1}{2} \log_{2} \left( 6p_{N}-2 \right)  
\end{equation}
now we have
\begin{equation}\label{eq:oddpowanyk}
	k_{j} = \left\lfloor \frac{1}{2} \log_{2} \left( \frac{6p_{N}-2}{6i_{opow}-1} \right) + \frac{1}{2} \right\rfloor
\end{equation}
and integer solutions for
\[ \frac{1}{2} \log_{2} \left( \frac{6p_{N}-2}{6i_{opow}-1} \right) = \frac{1}{2}(2f-1) \]
$f \in \mathbb{N}$.
\begin{equation}\label{eq:oddpowanyki}
	i_{opow} = \left\lfloor \frac{1}{6} \left( \frac{6p_{N}-2}{2^{-2f+1}} + 1 \right) \right\rfloor = \frac{1}{6} \left( \frac{6p_{N}-2}{2^{-2f+1}} + 1 \right) - R_{o}(p_{N},f)
\end{equation}
$0 \leq R_{o}(p_{N},f) < 1$. The maximum $i_{opow, max}$ for $f = 1$
\begin{equation}\label{eq:oddpowanykimax}
	i_{opow,max} = \left\lfloor \frac{1}{6} \left( \frac{6p_{N}-2}{2} + 1 \right) \right\rfloor = \left\lfloor \frac{1}{6} \left( 3p_{N} - 1 + 1 \right) \right\rfloor = \left\lfloor \frac{1}{2} p_{N} \right\rfloor
\end{equation}
So
\begin{equation} i_{opow,max} = \left\{
	\begin{array}{l@{\quad \quad}l}
	\frac{1}{2}p_{N} & \textnormal{if $p_{N}$ is even} \\
	\frac{1}{2}p_{N} - \frac{1}{2} & \textnormal{if $p_{N}$ is odd}
	\end{array}
\right. \label{eq:iopowmaxcases} \end{equation} 
\subsection{Even powers: qualities for fixed $N$}
\label{ss:evenqual}
For consideration of even powers we get for (\ref{eq:powrel}) with
\[ n_{2,i} = 1, x_{i} = 2k_{N} \]
\[ n_{2,j} = 6i_{epow} + 1, x_{j} = 2k_{j} \]
$k_{j} \in \mathbb{N}$.
\begin{equation}\label{eq:evpowk}
	k_{j} = \left\lfloor k_{N} + \frac{1}{2} \log_{2} \left( \frac{1}{6i_{epow}+1} \right) \right\rfloor
\end{equation}
Be $k_{N} \in \mathbb{N}: k_{N} > 1$. Equation (\ref{eq:evpowk}) has integer solutions for
\[ \log_{2}\left( \frac{1}{6i_{epow}+1} \right) = -2f \]
$f \in \mathbb{N}$.
\[ 2^{-2f} = \frac{1}{6i_{epow}+1} \]
\begin{equation}\label{eq:evpowi}
	i_{epow} = \left\lfloor \frac{1}{6} \left( \frac{1}{2^{-2f}} - 1 \right) \right\rfloor = \frac{1}{6} \left( \frac{1}{2^{-2f}} - 1 \right) - \frac{1}{2}
\end{equation}
and the maximum $i_{epow,max}$ for $f=k_{N}-1$
\begin{equation}\label{eq:evpowimax}
	i_{epow,max} = \left\lfloor \frac{1}{6} \left( \frac{1}{2^{-2k_{N}+2}} - 1 \right) \right\rfloor = \frac{1}{6} \left( \frac{1}{2^{-2k_{N}+2}} - 1 \right) - \frac{1}{2}
\end{equation}
At second we take (\ref{eq:anyN}). So
\begin{equation}\label{eq:evpowanyk}
	k_{j} = \left\lfloor \frac{1}{2} \log_{2} \left( \frac{6p_{N}-2}{6i_{epow}+1} \right) \right\rfloor
\end{equation}
and integer solutions for
\[ \log_{2} \left( \frac{6p_{N}-2}{6i_{epow}+1} \right) = 2f \]
$f \in \mathbb{N}$.
\begin{equation}\label{eq:epowanyki}
	i_{opow} = \left\lfloor \frac{1}{6} \left( \frac{6p_{N}-2}{2^{2f}} - 1 \right) \right\rfloor = \frac{1}{6} \left( \frac{6p_{N}-2}{2^{2f}} - 1 \right) - R_{e}(p_{N},f)
\end{equation}
$0 \leq R_{e}(p_{N},f) < 1$. The maximum $i_{epow, max}$ for $f = 1$
\begin{equation}\label{eq:epowanykimax}
	i_{epow,max} = \left\lfloor \frac{1}{6} \left( \frac{6p_{N}-2}{2^{2}} - 1 \right) \right\rfloor = \left\lfloor \frac{1}{6} \left( \frac{3p_{N}-3}{2} \right) \right\rfloor = \left\lfloor \frac{p_{N}-1}{2^{2}} \right\rfloor
\end{equation}
So
\begin{equation} i_{epow,max} = \left\{
	\begin{array}{l@{\quad \quad}l}
	\frac{p_{N}-1}{2^{2}} - \frac{1}{2^{2}} & \textnormal{if $p_{N} = 2,6,10,... = 4s-2$} \\
	\frac{p_{N}-1}{2^{2}} - \frac{2}{2^{2}} & \textnormal{if $p_{N} = 3,7,11,... = 4s-1$} \\
	\frac{p_{N}-1}{2^{2}} - \frac{3}{2^{2}} & \textnormal{if $p_{N} = 4,8,12,... = 4s$} \\
	\frac{p_{N}-1}{2^{2}} & \textnormal{if $p_{N} = 5,9,13,... = 4s+1$}
	\end{array}
\right. \label{eq:iepowmaxcases} \end{equation} 
$s \in \mathbb{N}$.
\subsection{Total of generated odd numbers}
\label{ss:totnum}
Now we show the last important quality of our number sequence before we will come to the final proof step. In section \ref{ss:nuitfor} we saw that every number, which is generated, is only by one possible combination for $(n_{2},x)$ generated. This in an important result, because we want to decree the count of odd numbers, which are generated until a fixed $N$ by pass through the number rows $n_{2}$ ordered by size. See also tables \ref{tab:nopoweven} and \ref{tab:nopoodd}.

We start with the subset $n_{2,opow} = 6i_{opow} - 1$ with odd powers. Under use of the results of last section, the total of odd numbers by this subset is given by
\begin{multline*} T_{o} = k_{N} \left( \frac{1}{6} \left( 2^{1}+1 \right) - \frac{1}{2} \right) + \\
+ \sum_{i=1}^{k_{N}-1} \left( k_{N} - i \right) \left( \left( \frac{1}{6} \left( 2^{2(i+1)-1} + 1 \right) - \frac{1}{2} \right) - \left( \frac{1}{6} \left( 2^{2i-1} + 1 \right) - \frac{1}{2} \right) \right) \end{multline*}
\[ = \sum_{i=1}^{k_{N}-1} \left( k_{N}-i \right) \left( \frac{1}{6} \left( 2^{2(i+1)-1} - 2^{2i-1} \right) \right) \]
Since (geometric sum formula)
\begin{equation}\label{eq:sum2}
	\sum_{i=a}^{b} 2^{2i} = \frac{1}{3}4^{b+1} - \frac{1}{3} 4^{a}
\end{equation}
and
\begin{equation}\label{eq:sumi2}
	\sum_{i=a}^{b} i2^{2i} = \frac{1}{3} 4^{b+1} \left( b+1 \right) - \frac{4}{9} 4^{b+1} - \frac{1}{3} 4^{a} a + \frac{4}{9} 4^{a}
\end{equation}	
it follows
\begin{equation}\label{eq:totodd}
	T_{o} = \frac{1}{9} 4^{k_{N}} - \frac{1}{3} k_{N} - \frac{1}{9}
\end{equation}
Here we used the results for $N = \frac{2^{2k_{N}}-1}{3}$, $k_{N} \in \mathbb{N}: k_{N} > 1$. In this context, the proof for this special case also includes the proof for $N=2p_{N}-1$, $p_{N} \in \mathbb{N}:p_{N} > 1$, since it is clearly that in the second case, we have the same rows and no other qualities for the total.

For the subset $n_{2,epow} = 6i_{epow} + 1$ with even powers it follows 
\begin{multline*} T_{e} =  \left( k_{N} - 1 \right) \left( \frac{1}{6} \left( 2^{2} - 1 \right) - \frac{1}{2} \right) + \\ 
+ \sum_{i=2}^{k_{N}-1} \left( k_{N}-i \right) \left( \left( \frac{1}{6} \left(2^{2i} - 1 \right) - \frac{1}{2} \right) - \left( \frac{1}{6} \left( 2^{2(i-1)} - 1 \right) - \frac{1}{2} \right) \right) \end{multline*}
\[ = \sum_{i=2}^{k_{N}-1} \left( k_{N} - i \right) \left( \frac{1}{6} \left( 2^{2i} - 2^{2(i-1)} \right) \right) \]
\begin{equation}\label{eq:toteven}
	T_{e} = \frac{1}{18} 4^{k_{N}} - \frac{2}{3} k_{N} + \frac{4}{9}
\end{equation}
The total of odd numbers in range $[1,N]$ which are generated by $n_{2}$ rows is given by
\begin{equation}\label{eq:totnumb}
	T = \left( k_{N} - 1 \right) + 1 + T_{o} + T_{e} = \frac{1}{6} 4^{k_{N}} + \frac{1}{3}
\end{equation}
Here, we have to considered the exception from section \ref{ss:it_opch} and the second exception from section \ref{ss:nuitfor}. We can show that this is equal to the total of all odd numbers in range $[1,N]$
\[ \frac{1}{6} 4^{k_{N}} + \frac{1}{3} \stackrel{!}{=} \left( \frac{2^{2k_{N}} - 1 + 3 }{3} \right) \frac{1}{2} =\frac{2^{2k_{N}-1}+1}{3} \]
\[ \frac{1}{2} 4^{k_{N}} + 1 \stackrel{!}{=} 2^{2k_{N}-1} + 1 \]
\[ \frac{1}{2} 2^{2k_{N}} \stackrel{!}{=} 2^{2k_{N}-1} \]
\[ 2^{2k_{N}-1} = 2^{2k_{N}-1} \]
$\forall k_{N} \in \mathbb{R}$.
\section{Recurrence of valid $N$ area}
\label{s:itofN}
Now we will make the final proof step under using the results of previous sections. Let's start with the following assumption
\begin{ass}\label{ass:finalass}
The trueness of the conjecture in a finite range $[1,N_{0}]$ was already shown. For example by explicit execution of generation rule for Collatz sequences. 
\end{ass}
So we will do this. For example, we want to start with the range $[1,19]$. To show the assumption, we generated the sequences, starting by the odd numbers within the range. See table \ref{tab:assumpt}.
\begin{table}
\centering
\begin{tabular}{|c|c|}
\hline $4 \rightarrow 2 \rightarrow$ \textbf{1} & $1$ \\ 
\hline \textbf{3} $\rightarrow 10 \rightarrow$ \textbf{5} $\rightarrow 16 \rightarrow 8 \rightarrow 4 \rightarrow \cdots$ & $3,5$ \\
\hline \textbf{7} $\rightarrow 22 \rightarrow$ \textbf{11} $\rightarrow 34 \rightarrow$ \textbf{17} $\rightarrow 52 \rightarrow 26 \rightarrow$ \textbf{13} $\rightarrow 40 \rightarrow 20 \rightarrow 10 \rightarrow \cdots$ & $7,11,13,17$ \\
\hline \textbf{9} $\rightarrow 28 \rightarrow 14 \rightarrow 7 \rightarrow \cdots$ & $9$ \\
\hline \textbf{15} $\rightarrow 46 \rightarrow$ \textbf{23} $\rightarrow 70 \rightarrow 35 \rightarrow 106 \rightarrow 53 \rightarrow 160 \rightarrow 80 \rightarrow 40 \rightarrow \cdots$ & $15$ \\
\hline \textbf{19} $\rightarrow 58 \rightarrow$ \textbf{29} $\rightarrow 88 \rightarrow 44 \rightarrow 22 \rightarrow 11 \rightarrow \cdots$ & $19$ \\
\hline 
\end{tabular}\caption{Collatz sequences started by the bold marked odd numbers off example range $[1,19]$. Also bold marked, all $6i-1$ numbers of new range $[1,29]$.}\label{tab:assumpt}\end{table}

Now we study the three subsets from above. At first, look at subset $n_{2,opow} = 6i_{opow} - 1$. It was shown that, for generating all numbers in the range $[1,N_{0}]$, we need all numbers of this set type until $i_{opow,max,0}$, so
\begin{equation}\label{eq:oddmax}
	n_{opow,max,0} = 6i_{opow,max,0} - 1 = 6 \left\lfloor \frac{1}{2} p_{N_{0}} \right\rfloor - 1 
\end{equation}
\begin{equation} n_{opow,max,0} = \left\{
	\begin{array}{l@{\quad \quad}l}
	6\frac{1}{2}p_{N_{0}}-1 = 3p_{N_{0}}-1 & \textnormal{if $p_{N_{0}}$ is even} \\
	6\left(\frac{1}{2}p_{N_{0}}-\frac{1}{2}\right) = 3p_{N_{0}}-4 & \textnormal{if $p_{N_{0}}$ is odd}
	\end{array}
\right. \label{eq:nopowmax0} \end{equation} 
We can show, that this new range from lowest point to $n_{opow,max,0}$ is always higher than the start number $N_{0}$.
\[ \Delta n_{0} = n_{opow,max,0} - N_{0} = \left\{
	\begin{array}{l@{\quad \quad}l}
	3p_{N_{0}}-1 - 2p_{N_{0}} + 1 = p_{N_{0}} & \textnormal{if $p_{N_{0}}$ is even} \\
	3p_{N_{0}}-4 - 2p_{N_{0}} + 1 = p_{N_{0}} - 3 & \textnormal{if $p_{N_{0}}$ is odd}
	\end{array}
\right. \] 
For our example this means the following. For the range $[1,19]$, $p_{N_{0}} = 10$
\[ n_{opow,max,0} = 3p_{N_{0}}-1 = 29 \]
Numbers which are element of this subset are
\[ 5, 11, 17, 23, 29 \]
In table \ref{tab:assumpt} we can see that all of this numbers have to pass as direct precursor of the final range numbers on the way from number one to them. So we have shown, that, for generating all numbers from $[1,N_{0}]$, before, we have to go over all numbers of the subset $6i_{opow,0} - 1$ in the higher finite range $[1,n_{opow,max,0}]$. Our new range for, only numbers of this subset, is so
\begin{equation}\label{eq:oddRange}
	N_{o,1} = n_{opow,max,0}
\end{equation}
Now we will look at subset $n_{2,epow,0} = 6i_{epow,0} + 1$. Here we have
\begin{equation}\label{eq:evmax}
	n_{epow,max,0} = 6i_{epow,max,0} + 1 =  6 \left\lfloor \frac{p_{N_{0}}-1}{2^{2}} \right\rfloor + 1
\end{equation}
\begin{equation} n_{epow,max,0} = \left\{
	\begin{array}{l@{\quad \quad}l}
	6\left( \frac{p_{N_{0}}-1-1}{2^{2}} \right) + 1 = \frac{3p_{N_{0}}-4}{2} & p_{N_{0}}=4s_{0}-2 \\
	6\left( \frac{p_{N_{0}}-1-2}{2^{2}} \right) + 1 = \frac{3p_{N_{0}}-7}{2} & p_{N_{0}}=4s_{0}-1 \\
	6\left( \frac{p_{N_{0}}-1-3}{2^{2}} \right) + 1 = \frac{3p_{N_{0}}-10}{2} & p_{N_{0}}=4s_{0} \\
	6\left( \frac{p_{N_{0}}-1}{2^{2}} \right) + 1 = \frac{3p_{N_{0}}-1}{2} & p_{N_{0}}=4s_{0}+1 \\	
	\end{array}
\right. \label{eq:epowmax0cases} \end{equation} 
$s_{0} \in \mathbb{N}$. But now, with
\[ \Delta n_{0} = n_{epow,max,0} - N_{0} = \left\{
	\begin{array}{l@{\quad \quad}l}
	\frac{3p_{N_{0}}-4}{2} - 2p_{N_{0}} + 1 = \frac{-p_{N_{0}}-2}{2} & p_{N_{0}}=4s_{0}-2 \\
	\frac{3p_{N_{0}}-7}{2} - 2p_{N_{0}} + 1 = \frac{-p_{N_{0}}-5}{2} & p_{N_{0}}=4s_{0}-1 \\
	\frac{3p_{N_{0}}-10}{2} - 2p_{N_{0}} + 1 = \frac{-p_{N_{0}}-8}{2} & p_{N_{0}}=4s_{0} \\
	\frac{3p_{N_{0}}-1}{2} - 2p_{N_{0}} + 1 = \frac{-p_{N_{0}}+1}{2} & p_{N_{0}}=4s_{0}+1 \\
\end{array}
\right. \] 
it always follows $\Delta n_{0} < 0$. We can solve this problem, when we use an other angle of view for our numbers. We see that all numbers in range $[1,n_{epow,max,0}]$ from this subset are generated. Also we know, from assumption, that all numbers of this subset are generated in $[n_{epow,max,0},N_{0}]$. So we know, that all numbers of this subset can be considered as a new $n_{epow,max,1} = N_{0}$ for a bigger range $N_{e1}$. So we have
\begin{equation}\label{eq:pe1an}
	2p_{N_{0}} - 1 = N_{0} = 6 \left\lfloor \frac{p_{N_{1}} - 1}{2^{2}} \right\rfloor + 1
\end{equation}
This equation has to satisfy two needs. On the one hand, $p_{N_{1}}$ has to be an integer number and otherwise, $\frac{p_{N_{1}}-1}{2^{2}}$ has to be an integer number, too.
\[ 8p_{N_{0}} - 4 = 6p_{N_{1}} - 2 \]
\begin{equation}\label{eq:pe1}
	p_{N_{1}} = \frac{4p_{N_{0}} - 1}{3}
\end{equation}
\begin{equation} p_{N_{1}} = \left\{
	\begin{array}{l@{\quad \quad}l}
	\frac{4p_{N_{0}}-1-1}{3} = 4s-2 & p_{N_{0}}=2,5,8,... = 3s_{0} - 1 \\
	\frac{4p_{N_{0}}-1-2}{3} = 4s-1 & p_{N_{0}}=3,6,9,... = 3s_{0} \\
	\frac{4p_{N_{0}}-1}{3} = 4s+1& p_{N_{0}}=4,7,10,... = 3s_{0} + 1  
\end{array}
\right. \label{eq:pn1cases} \end{equation} 
$s_{0} \in \mathbb{N}$. With equation (\ref{eq:pe1an})
\[ \frac{p_{N_{1}}-1}{2^{2}} = \left\{
	\begin{array}{l@{\quad \quad}l}
	\frac{4s-3}{2^{2}} - \frac{1}{4} & \Rightarrow p_{N_{0}} = 1,4,7,... = 3s-2 \\
	\frac{4s-2}{2^{2}} - \frac{2}{4} & \Rightarrow p_{N_{0}} = 1,4,7,... = 3s-2 \\
	\frac{4s}{2^{2}} & \Rightarrow p_{N_{0}} = 4,7,10,... = 3s+1   
\end{array}
\right. \] 
To solve equation (\ref{eq:pe1an}) correct, we need
\[ \frac{p_{N_{1}}-1}{2^{2}} = \left\{
	\begin{array}{l@{\quad \quad}l}
	\frac{4s-3}{2^{2}} - \frac{1}{4} + \frac{1}{3} & p_{N_{0}} = 3s - 1 \\
	\frac{4s-2}{2^{2}} - \frac{2}{4} + \frac{2}{3} & p_{N_{0}} = 3s \\
	\frac{4s}{2^{2}} & p_{N_{0}} = 3s+1   
\end{array}
\right. \] 
We get for the range $N_{e,1} = 2p_{N_{1}}-1$
\begin{equation} N_{e,1} = \left\{
	\begin{array}{l@{\quad \quad}l}
	8s_{0}-5 = \frac{8p_{N_{0}}-7}{3} = \frac{4N_{0}-3}{3} & p_{N_{0}}=2,5,8,... = 3s_{0} - 1 \\
	8s_{0}-3 = \frac{8p_{N_{0}}-9}{3} = \frac{4N_{0}-5}{3} & p_{N_{0}}=3,6,9,... = 3s_{0} \\
	8s_{0}+1 = \frac{8p_{N_{0}}-5}{3} = \frac{4N_{0}-1}{3} & p_{N_{0}}=4,7,10,... = 3s_{0} + 1  
\end{array}
\right. \label{eq:Ne1} \end{equation}
Now we have a look at our example $p_{N_{0}} = 10$. For the subset $n_{2,epow,0} = 6i_{epow,0} + 1$ we have
\[ N_{e,1} = \frac{8 \cdot 10 - 5}{3} = 25 \]
Numbers which are element of this subset are
\[ 7, 13, 19, 25 \]
The numbers $7,13$ are elements of $[1,N_{0}]=[1,19]$ and $19,25$ are elements of the new range $[1,N_{1}]=[1,25]$. See also tables \ref{tab:nopoweven}, \ref{tab:nopoodd} and \ref{tab:assumpt}.
To find our final new range we have to compare the two calculated ranges $N_{0,1}$ and $N_{e,1}$ and choose the smaller of them. For this let be
\[ \Delta N_{oe,1} = N_{o,1} - N_{e,1} = \left( 3p_{N_{0}} - A \right) - \left( \frac{8p_{N_{0}}-B}{3} \right) \]
$A := \{1,4\}$, $B := \{ 7,9,5 \}$. 
\[ = \frac{9p_{N_{0}}-3A-8p_{N_{0}}+B}{3} \]
\begin{equation} \Delta N_{oe,1} = \frac{p_{N_{0}}+B-3A}{3} \label{eq:deltanoe1} \end{equation}
for $A=4$ and 
\begin{equation} \Delta N_{oe,1} = \left\{
	\begin{array}{l@{\quad \quad}l}
	\frac{p_{N_{0}}-5}{3} \rightarrow \Delta N_{oe,1} = 0 \rightarrow N_{o,1} = N_{e,1} & B=7, p_{N_{0}}=5 \\
	\frac{p_{N_{0}}-3}{3} \rightarrow \Delta N_{oe,1} < 0 \rightarrow N_{o,1} < N_{e,1} & B=9, p_{N_{0}}=3 \\
	\frac{p_{N_{0}}-7}{3} \rightarrow \Delta N_{oe,1} = 0 \rightarrow N_{o,1} = N_{e,1} & B=5, p_{N_{0}}=7  
\end{array}
\right. \label{eq:deltaone1cases}\end{equation} 
else
\begin{equation} \Delta N_{oe,1} > 0 \Rightarrow N_{0,1} > N_{e,1} \label{eq:deltaone1pos} \end{equation}
We have to use the smaller one of this two ranges for the new recurrence range $N_{1}$. For the cases above we have
\begin{equation} N_{1} = \left\{
	\begin{array}{l@{\quad \quad}l}
	N_{o,1} = N_{e,1} & B=7, p_{N_{0}}=5 \\
	N_{o,1} & B=9, p_{N_{0}}=3 \\
	N_{o,1} = N_{e,1} & B=5, p_{N_{0}}=7  
\end{array}
\right. \label{eq:n1cases}\end{equation} 	
Else it follows
\begin{equation}\label{eq:N}
	N_{1} = N_{e,1}
\end{equation}
We know that all numbers of the subsets $6i_{opow}-1$ and $6i_{epow}+1$ are elements of this new range.

Finally, what is about the third subset? The numbers which are integer multiple numbers of three. Always these numbers are all generated by the other two subsets for our new range $N_{1}$, which we have already shown in section \ref{ss:totnum}. From this it follows
\begin{center}
	We have a new range $N_{1} > N_{0}$ in which are all odd numbers can be reached by number one. 
\end{center}
To show that we can do this recurrence any often, the range $N_{i+1}$ has to be greater than $N_{i}$ for all recurrence steps $i \in \mathbb{N}$. This is nearly always true, since from equation (\ref{eq:nopowmax0}) 
\[ N_{o,i+1} = \left\{
	\begin{array}{l@{\quad \quad}l}
	3p_{N_{i}}-1 = \frac{3}{2}N_{i} + \frac{1}{2} & \textnormal{if $p_{N_{i}}$ is even} \\
	3p_{N_{i}}-4 = \frac{3}{2}N_{i} - \frac{5}{2} & \textnormal{if $p_{N_{i}}$ is odd}
	\end{array}
\right. \] 
We see that always $N_{o,i+1} > N_{i}$ apart from the case $p_{N_{i}} = 3$, in which we have
\[ N_{i+1} = N_{i} \]
And from equation (\ref{eq:Ne1}) we have
\[ N_{e,i+1} = \frac{4N_{i}-C}{3} \]
$C := \{ 3, 5, 1 \}$.
\[ = \frac{4}{3}N_{i} - C \]
which gives us the exceptions $p_{N_{i}} = 2$ and $p_{N_{i}} = 3$. So we have only the cases $p_{N_{i}}=2$ and $p_{N_{i}}=3$ in which we can't do the correct recurrence. But this problem is easy solved. If we used for our assumption $p_{N_{0}} > 3$, for which we can show, for example by building the concrete number sequence, the validity of Collatz conjecture, we can do our iterations for all numbers which are greater. 

So, we see that we can do this recurrence endless. $\Rightarrow$ The conjecture is true for all odd numbers. $\Rightarrow$ The conjecture is true for all even numbers. $\Rightarrow$ The conjecture is true for all $N \in \mathbb{N}$.
\begin{flushright} $\Box$ \end{flushright}
\section{Closing notes.}
\label{s:closnot}
Closing at this point I want to make some notes on my proof and the content of two publications which I have found. In \cite{Cpunp} and \cite{ComSS}, Craig Alan Feinstein gives a proof why the Collatz conjecture be unprovable. I deliberate about this proof and my solution a lot of time. I think it's important to point on this example to show that there can be issues of mathematical structures which can be forgot when we look at it. It's a good example that we should never be too confident by using our mathematical tool kit.
\nocite{*}
\bibliography{literature}
\bibliographystyle{aomalpha}

\end{document}